\documentclass[12pt]{amsart}
\usepackage{amssymb}
\usepackage[all]{xy}
\usepackage[pdftitle={finite Heisenberg groups}]{hyperref}

 \rm

\newcommand{\x}{\times}

\newcommand{\abs}[1]{\left\lvert#1\right\rvert}

\newcommand{\st}{\:|\:}

\newcommand{\C}{{\mathbb{C}}}

\newcommand{\Z}{{\mathbb{Z}}}

\renewcommand{\phi}{\varphi}
\newcommand{\p}{\rho}

\renewcommand{\H}{{\mathcal{H}}}
\newcommand{\BH}{{\mathcal{B}}(\H)}
\newcommand{\A}{{\mathcal{A}}}

\newcommand{\End}{{\mathrm{End}}}

\theoremstyle{plain}
\newtheorem{thm}{Theorem}[section]
\newtheorem{lem}[thm]{Lemma}

\theoremstyle{definition}
\newtheorem{defn}[thm]{Definition}

\theoremstyle{remark}

\title{Inductive algebras for finite Heisenberg groups}

\author{Amritanshu Prasad}
\address{The Institute of Mathematical Sciences, CIT campus Taramani, Chennai 600113.}

\author{M.~K.~Vemuri}
\address{Chennai Mathematical Institute, Plot H1, SIPCOT IT Park, Padur~PO, Siruseri 603103.}


\begin{document}

\begin{abstract}
A characterization of the maximal abelian sub-algebras of matrix
algebras that are normalized by the canonical representation of
a finite Heisenberg group is given.  Examples are constructed using a classification result for finite Heisenberg groups.
\end{abstract}

\maketitle


\section{Introduction}\label{S:intro}

Let $G$ be a separable locally compact group and $\p$ an irreducible
unitary representation of $G$ on a separable Hilbert space $\H$.  Let
$\BH$ denote the algebra of bounded operators on $\H$.  An {\em
inductive algebra} is a weakly closed abelian sub-algebra $\A$ of $\BH$
that is normalized by $\p(G)$, i.e., $\p(g) \A \p(g)^{-1} = \A$ for each
$g \in G$.  If we wish to emphasize the
dependence on $\p$, we will use the term $\p$-inductive algebra.
A {\em maximal inductive algebra} is a maximal element of the set of
inductive algebras, partially ordered by inclusion.

The identification of inductive algebras can shed light on the
possible realizations of $\H$ as a space of sections of a homogeneous
vector bundle (see, e.g., \cite{Stegel-I,Stegel-II,Vemuri-SL2R.,Vemuri-RCR}).  For self-adjoint maximal inductive algebras there is a
precise result known as Mackey's Imprimitivity Theorem, as explained
in the introduction to \cite{Stegel-I}.

Inductive algebras were introduced in \cite{Vemuri-RCR}, where the
maximal inductive algebras for the real Heisenberg groups were
classified.
In this paper, we classify the maximal inductive algebras for finite Heisenberg groups.
All these algebras are found to be self-adjoint.
By Mackey's imprimitivity theorem, they arise from induction.
Therefore, we get a classification of the realizations of the canonical representation as induced representations.
Our proof of the classification of maximal inductive algebras in the finite case is very direct and reveals the simple algebraic idea behind such a result.

In Section \ref{S:summary}, we recall the representation theory of
a general Heisenberg group, and formulate our main result.  Our definition
of Heisenberg group is taken from
\cite{Theta-III}.  The main result is proved in
Section~\ref{S:iafhg}.
In Section~\ref{sec:examples}, we give a classification of finite Heisenberg
groups and use it to construct a large class of maximal inductive algebras
for their canonical representations.

\section{Formulation of the Main Theorem}\label{S:summary}

Let $K$ be a locally compact abelian group, and 
$U(1) =\{ z \in \C \: :\: \abs{z} = 1 \}$.
Consider a locally compact group $G$ which lies in a central extension:
\begin{equation}\label{hg}
\xymatrix{
1 \ar[r] & U(1) \ar[r] & G \ar[r]^\pi & K \ar[r] & 0
}
\end{equation}
Assume that $\pi$ admits a continuous section $s$.  Thus $s: K \to G$
and $\pi \circ s = {\mathrm{Id}}_K$.  Define $e: K \x K \to U(1)$ by
$e(k,l) = s(k)s(l)s(k)^{-1}s(l)^{-1}$.  Then $e$ is well
defined and independent of the choice of the section $s$.  Furthermore,
$e$ is an alternating bicharacter, i.e., it is a homomorphism in each
variable separately, and $e(k,k) = 1$ for each $k\in K$.  
Finally, $e$ determines $G$ up to unique isomorphism of central extensions.
Let $\widehat{K}$ denote
the Pontryagin dual of $K$.  Then a homomorphism $e^\flat: K \to \widehat{K}$
is obtained by setting
$$
(e^\flat(l))(k) = e(k,l).
$$

\begin{defn}
If $e^\flat$ is an isomorphism, then $G$ is called a {\em Heisenberg group}.
\end{defn}

Let $G$ be a Heisenberg group as in (\ref{hg}).  Then, $Z(G)=U(1)$.
We will assume that $K$ is finite.
The resulting Heisenberg groups are customarily called \emph{finite Heisenberg groups}.

\begin{defn}
Let $M \subseteq K$ be a closed subgroup.  If $e|_{M \x M} \equiv 1$ then
$M$ is called an {\em isotropic subgroup}.  A {\em maximal isotropic
subgroup} is one which is a maximal element of the set of isotropic
subgroups, partially ordered by inclusion.
\end{defn}

Note that there is a one to one correspondence between isotropic
subgroups of $K$ and abelian subgroups of $G$ containing $Z(G)$ given by
$M \mapsto \pi^{-1}(M)$.  Consequently, there is a one to one
correspondence between maximal isotropic subgroups of $K$ and
maximal abelian subgroups of $G$.

If $M \subseteq K$ is isotropic then we have an exact sequence of
abelian groups
\begin{equation*}
\xymatrix{
1 \ar[r] & U(1) \ar[r] & \pi^{-1}(M) \ar[r]^\pi & M \ar[r] & 0
}
\end{equation*}
The exactness of the dual sequence
\begin{equation*}
\xymatrix{
0 & \ar[l] \Z & \ar[l] \widehat{\pi^{-1}(M)} & \ar[l] \widehat{M} & \ar[l] 0
}
\end{equation*}
implies that the character $\mathrm{id}_{Z(G)}: Z(G) \to U(1)$ extends to a
character of $\pi^{-1}(M)$.

Let $M \subseteq K$ be a maximal isotropic subgroup, and let
$\chi \in \widehat{\pi^{-1}(M)}$ be such that $\chi |_{Z(G)} = \mathrm{id}_{Z(G)}$.
Then the representation $\p = \mathrm{Ind}_{\pi^{-1}(M)}^G \chi$ is irreducible, and
called the {\em canonical representation}.  The following is a
famous theorem of Stone, von Neumann and Mackey.

\begin{thm}
If $\p'$ is an irreducible representation of $G$ on a vector
space $V$ such that $\p'|_{Z(G)} = (\mathrm{id}_{Z(G)}) \cdot I$ then
$\p' \cong \p$.
\end{thm}

A proof of this theorem may be found in \cite{Theta-III}.
We may think of
$\p$ as acting on the space
$$
V = \{ f:G \to \C \st f(hg')=\chi(h)f(g') \text{ for all } h \in \pi^{-1}(M),\; g' \in G \}
$$
by
$ (\p(g)f)(g')=f(g'g)$.  If $u, v: G \to \C$ are functions, define
$m_u(v) = uv$.  Thus $m_u$ is the operator of multiplication by $u$.
Let
$$
\A(M) = \{ m_u \st u(hg') = u(g') \text{ for all } h \in \pi^{-1}(M),\; g' \in G \}.
$$
Then $\A(M) \subseteq \End(V)$ is a maximal abelian subalgebra.  Furthermore,
it is inductive.  So it is a maximal inductive algebra.  We note that
$\A(M)$ is self adjoint.  Our main result is the following.

\begin{thm}
  \label{thm:main}
  Every maximal $\p$-inductive algebra is of the form $\A(M)$ for a maximal
  isotropic subgroup $M$ of $K$.  In particular, all maximal $\p$-inductive
  algebras are maximal abelian subalgebras and are self-adjoint.
\end{thm}

\section{The Proof}\label{S:iafhg}

For $k \in K$, define $\kappa(k): \End(V) \to \End(V)$ by
by $\kappa(k)T = \p(s(k)) T \p(s(k))^{-1}$.  Then $\kappa(k)$ is
independent of the choice of section $s$, and 
$\kappa: K \to \mathrm{GL}(\End(V))$
is a representation of the finite abelian group $K$.
For each $l\in K$, let 
\begin{equation*}
  \End(V)_l=\{T\in\End(V) \;|\: \kappa(k)(T)=e^\flat(l)(k)T\}.
\end{equation*}

\begin{lem}
  \label{lem:eigenvectors}
  $\End(V)$ has a multiplicity-free decomposition into $K$-invariant subspaces
  \begin{equation*}
    \End(V)=\bigoplus_{l\in K} \End(V)_l.
  \end{equation*}
  Moreover, the subspace $\End(V)_l$ is spanned by $\p(s(l))$.
\end{lem}
\begin{proof}
  Since
  \begin{eqnarray*}
    \kappa(k)(\p(s(l))) & = & \p(s(k))\p(s(l))\p(s(k))^{-1}\\
    & = & \p(s(k)s(l)s(k)^{-1} s(l)^{-1}) \p(s(l))\\
    & = & e^\flat(l)(k)\p(s(l)),
  \end{eqnarray*}
  $\p(s(l))\in \End(V)_l$.
  The vectors $\p(s(l))$, $l\in K$ are linearly independent, as they lie in distinct eigenspaces for $K$.
  Since $\dim(\End(V))=|K|$, they must form a basis of $\End(V)$, and the result follows.
\end{proof}
To prove Theorem~\ref{thm:main}, note that every maximal $\p$-inductive algebra $\A$ is an invariant subspace of $\End(V)$ under the action of $K$ by $\kappa$.
By Lemma~\ref{lem:eigenvectors}, such a subspace is the direct sum of those $\End(V)_l$ which are contained in it.
Hence, $\A$ is determined by the set
\begin{equation*}
  M=\{l\in K\;|\: \p(s(l))\in \A\}.
\end{equation*}
The commutativity of $\A$ implies that for any $l,k\in M$, $e(l,k)=1$.
The maximality of $\A$ implies that $M$ is a maximal subset of $K$ with this property. 
In other words, $M$ is a maximal isotropic subgroup of $K$.

For each $g\in G$, there exists $z_g\in U(1)$ such that $g=z_gs(\pi(g))$. 
Now, for any $f\in V$,
\begin{eqnarray*}
  \p(s(l))f(g) & = & f(z_g\pi(g)s(l))\\
  & = & f(z_g (\pi(g)s(l)\pi(g)^{-1}s(l)^{-1})s(l)\pi(g))\\
  & = & e(s(l),\pi(g))\chi(s(l))f(g).
\end{eqnarray*}
Hence, $\p(s(l))=m_u$, where $u(g)=e(s(l),\pi(g))\chi(s(l))$.
Since $m_u$ preserves $V$, $u(hg)=u(g)$ for all $h\in \pi^{-1}(M)$ and $g\in G$.
Hence $\rho(s(l))\in \A(M)$ for each $l\in M$, and therefore, $\A\subset \A(M)$.
Since $e^\flat|_M$ descends to an isomorphism $K/M\to \hat M$,
\begin{equation*}
  \dim \A=|M|=|\hat M|=|K/M|=\dim \A(M).
\end{equation*}
Therefore, $\A=\A(M)$.
\section{Examples}
\label{sec:examples}
We shall construct examples of maximal isotropic subgroups, and hence of maximal inductive algebras, using a classification the finite Heisenberg groups up to isomorphism.

Since the pair $(K,e)$ as in Section~\ref{S:summary} determines the Heisenberg group $G$ up to isomorphism, it suffices to classify such pairs.
If $A$ is a locally compact abelian group denote its Pontryagin dual by $\hat A$.
\begin{thm}
  \label{thm:factorization}
  Suppose that $K$ is a finite abelian group and $e:K\times K\to U(1)$ is a non-degenerate alternating bicharacter.
  Then there exists a finite abelian group $A$ and an isomorphism $\phi:A\times \hat A\to K$ such that
  \begin{equation*}
    e(\phi(x,\chi),\phi(x',\chi'))=\chi'(x)\chi(x')^{-1} \text{ for all }x,x'\in A\text{ and } \chi,\chi'\in \hat A.
  \end{equation*}
\end{thm}
\begin{proof}
  Choose a basis $x_1,\ldots,x_n$ of $K$ such that the order of $x_i$ is $d_i$ with $d_1|d_2|\cdots|d_n$.
  Let $\zeta$ denote a fixed primitive $d_1$th root of unity.
  For each $i,j\in \{1,\ldots,n\}$, there exists a non-negative integer $q_{ij}$ such that
  \begin{equation*}
    e(x_i,x_j)=\zeta^{q_{ij}}.
  \end{equation*}
  Observe that $q_{ij}$ is determined modulo $d_1$ and is constrained to be divisible by $d_1/(d_i,d_j)$.
  Let $Q(e)$ denote the matrix $(q_{ij})$ associated to $e$ as above.
  
  Any automorphism $\alpha$ of $K$ is of the form
  \begin{equation*}
    \textstyle \alpha\big(\sum_{j=1}^n a_jx_j\big) = \sum_{i=1}^n \sum_{j=1}^n \alpha_{ij} a_j x_i
  \end{equation*}
  for a suitable matrix $(\alpha_{ij})$ of integers,
  where $\alpha_{ij}$ is determined modulo $d_i$ and is constrained to be divisible by $d_i/(d_i,d_j)$.
  If $e'(k,l)=e(\alpha(k),\alpha(l))$, then 
  \begin{equation*}
    Q(e')={}^t\alpha Q(e) \alpha.
  \end{equation*}
  In the above expression, $\alpha$ is identified with its matrix $(\alpha_{ij})$.
  In particular, $Q(e')$ can be obtained from $Q(e)$ by a sequence of column operations followed by the corresponding sequence of row operations.
  Following Birkhoff \cite{Birkhoff}, we use the operations
  \begin{enumerate}
  \item [$\beta_{i\sigma}$: ] multiplication of the $i$th row and $i$th column by an integer $\sigma$ such that $(\sigma,d_i)=1$.
  \item [$\pi_{ij}$: ] interchange of $i$th and $j$th rows and columns when $d_i=d_j$.
  \item [$\alpha_{ij\sigma}$: ] addition of $\sigma$ times $i$th row to the $j$th row, and $\sigma$ times the $i$th column to the $j$th column when $d_i/(d_i,d_j)$ divides $\sigma$.
  \end{enumerate}

  Since $e$ is alternating, $Q(e)$ can be taken to be skew-symmetric with $0$'s along the diagonal.
  Since $e$ is non-degenerate, there must be an integer in the first row of $Q(e)$ coprime to $d_1$.
  First, $\beta_{i\sigma}$ for appropriate $\sigma$ can be used to make this entry equal to $1$.
  If this entry is in the $j$th column, then the constraints on $q_{ij}$ force $d_j=d_1$.
  Now, $\pi_{2j}$ transforms the top-left corner of the matrix into $\left(\begin{smallmatrix} 0 & 1 \\ -1 & 0 \end{smallmatrix}\right)$.
  In particular $d_1=d_2$.
  Finally, $\alpha_{1j\sigma}$ and $\alpha_{2j\sigma}$ can be used to make the remaining entries of the first two rows and columns vanish.

  The form $e$ identifies the Pontryagin dual of the subgroup generated by the first (transformed) generator with the subgroup generated by the second one.
  After splitting off the subgroup spanned by these two generators, we are left with a finite abelian group with $n-2$ generators.
  Thus Theorem~\ref{thm:factorization} is proved by induction on $n$.
\end{proof}
Let $G$ be a finite Heisenberg group as in (\ref{hg}) and $\phi$ and $A$ be as in Theorem~\ref{thm:factorization}.
For any subgroup $B$ of $A$, let $B^\perp\subset \hat A$ denote the subgroup of characters which vanish on $B$.
Then $\phi(B\times B^\perp)$ is a maximal isotropic subgroup of $K$.
This class of maximal isotropic subgroups is likely to be very rich since the embedding problem for finite abelian groups involves wild phenomena \cite{Wild}.
\bibliographystyle{amsplain}
\bibliography{v10-iafhg}

\end{document}